\documentstyle[12pt,bezier,leqno]{article}
\begin{document}

\begin{center}
{\Large\sc On the Monodromy at Infinity }\\
{\Large\sc of a Polynomial Map, II}\\
$\ $ \\
{R. Garc\'\i a L\'opez\footnote{Supported by the DGICYT
}
and A. N\'emethi\footnote{Partially supported by OSU Seed Grant}}
\end{center}

\begin{footnotesize}
Universidad de Barcelona, Dept. de Algebra y Geometr\'{\i}a. Gran Via, 585,
08007 Barcelona, Spain.

The Ohio State University, 231 West 18th Avenue, Columbus OH 43210, USA.
\end{footnotesize}

\newcommand{\ea}{{\bf C}^{n+1}}
\newcommand{\cpx}{{\bf C}}
\newcommand{\zet}{{\bf Z}}
\newcommand{\pr}{{\bf C}[X_1,\ldots,X_{n+1}]}
\newcommand{\pe}{{\bf P}^{n+1}}
\newcommand{\ppe}{{\bf P}^n}
\newcommand{\ls}{j_{\ast}{\bf V}}
\newcommand{\cls}[1]{H^{#1}(\ppe, j_{\ast}{\bf V}_s)}
\newcommand{\mi}{T^{\infty}_{f}}
\newcommand{\xo}{X_{f}^{\infty}}
\newcommand{\ov}[1]{\overline{#1}}
\newcommand{\ti}[1]{\widetilde{#1}}
\newcommand{\xy}{X^{\infty}}
\newcommand{\vs}{\vspace{6mm}}
\newcommand{\vsss}{\vspace{4mm}}
\renewcommand{\labelenumi}{\alph{enumi})}
\renewcommand{\labelenumii}{\roman{enumii})}

\vs

\noindent {\bf \S 1. Introduction}
\vs

In the last years a lot of work has been
concentrated on the study
of the behaviour at infinity of polynomial maps
(see for example \cite{Lib}, \cite{Par}, \cite{Nem},
\cite{NZ}, \cite{GN}, among
others). This behaviour can be very complicated, therefore the main idea
was to find special classes of
polynomial maps which have, in some sense, nice properties at infinity.
In this paper, we completely determine the complex algebraic monodromy
at infinity for a special class of polynomial maps (which is
complicated enough to show the nature of the general problem).

Next, we give the precise definitions: Let $f:\ea \to \cpx$ be a map
given by a polynomial with complex coefficients (which will be
also denoted by $f$).
Then there exists a finite set $\Gamma \subset \cpx$ such that the map
\[
f\Bigm|_{ \ea - f^{-1}(\Gamma )}:\ea - f^{-1}(\Gamma ) \to
\cpx - \Gamma
\]
is a locally trivial ${\cal C}^\infty$-fibration (\cite{Pham}). We
denote by $\Gamma _f$ the smallest subset of the complex plane with
this property.
$\Gamma _f$ contains the set $\Sigma _f$ of critical values of $f$, but
in general it is bigger. Fix $t_0\in\cpx$ such that $\left| t_0 \right|>
\mbox{max}\{\left|t\right| : t\in\Gamma _f \}$. The complex algebraic monodromy
associated with the path $s\mapsto t_0 e^{2\pi i s}$, $s\in[0,1]$, is
denoted by
\[
(\mi)^{\ast}: H^{\ast}(f^{-1}(t_0),\cpx)\to H^{\ast}(f^{-1}(t_0),\cpx).
\]
This isomorphism is called the monodromy at infinity of $f$.
As we will see later, $(\mi)^{\ast}$ is a very delicate invariant of
$f$.

On studying topological properties of polynomial maps,
one usually imposes
some condition which insures the absence of vanishing cycles ``at infinity''
for a suitable compactification of the map $f$
(tameness, Malgrange condition,... cf. \cite{Par}). From this point of
view, a class of polynomial maps which looks natural to study is the
following:
\vs

\noindent {\it Definition:} A polynomial $f\in\pr$ will be called a
$(\ast )$-polynomial if it verifies the following condition:
\[
(\ast)
\left\{
\begin{array}{l}
\mbox{{\it For $t\in \cpx - \Sigma_f$, the closure in $\pe$ of the
affine}}\\
\mbox{{\it hypersurface $\{ f=t \}$ is non-singular.}}
\end{array} \right.
\]
The goal of this article is the computation of $(\mi)^{\ast}$ for
$(\ast)$-polynomials.\\
We will assume that $n\geq 2$. The case $n=1$ is completely clarified in
\cite{GN}.
\vs

If $d=$deg$(f)$ and $f=f_d+f_{d-1}+...$ is the decomposition of $f$ into
homogeneous components, condition $(\ast)$ is equivalent to
\[
\{ x\in \ea \mid \mbox{grad} f_d(x)=0, \ f_{d-1}(x)=0 \}=\{ 0 \},
\]
where grad denotes the gradient vector.
In the first part \cite{GN} of this sequence of papers, the following
results are given (besides others):
\begin{enumerate}
\item A $(\ast)$-polynomial $f$ satisfies $\Gamma _f = \Sigma _f$ and
any fiber of $f$ has the homotopy type of a bouquet of $n$-dimensional
spheres (cf.  \cite{Di.Top}). In particular,
the only interesting monodromy transformation is $(\mi)^n$, which in
the sequel will be denoted simply by $\mi$.
\item The hypersurface $\xy\subset \ppe$ given by $f_d=0$ has only
isolated singularities, and the monodromy at infinity (actually, the
whole topology at infinity) depends only on the hypersurface $\xy$.
\item The characteristic polynomial of $\mi$ is computable in terms of
the characteristic polynomials of the local monodromies of the
isolated singularities of $\xy$ (cf. Corollary 2).
\item On the other hand, the nilpotent part of $\xy$ cannot be
determined only from local data attached to the isolated singularities
of $\xy$, it depends essentially on the {\it position} of these
singular points.
\end{enumerate}

Part of the {\em global information} about the position of the singular points
of $\xy$ is already encoded in its Betti numbers. More subtle invariants
are hidden in the  complement $\ppe - \xy$ of
$\xy$, or in the cyclic coverings of $\ppe$ branched along $\xy$. For
algebraic surfaces, O. Zariski related this kind of invariants with the
defect (or superabundance)
of some linear systems, respectively with some Betti
numbers of cyclic coverings.
In the sequel we give the numerical invariants of $\xy$ which will
provide our description of $\mi$.

For $X$ a quasi-projective variety, denote by $b_q(X)$ (respectively,
$p_q(X)$) the dimension of $H^q(X,\cpx)$ (respectively, the dimension
of the $q$-th primitive cohomology of $X$). The numbers
$p_n(\xy)=b_{n-1}(\ppe -\xy)$ and $p_{n-1}(\xy)=b_{n}(\ppe - \xy)$ are
in general global invariants of $\xy$ (Here, if $n=2$, we define
$p_2(\xy)=b_2(\xy)-1)$). We define a map $h:\pi _1 (\ppe - \xy)\to
\zet /d\zet$ as follows: If $n>2$ then $h$ is just the Hurewicz map
(in fact, isomorphism):
\[
\pi _1 (\ppe - \xy) \to H_1(\ppe - \xy, \zet)=\zet / d\zet .
\]
If $n=2$, let $r$ denote the number of irreducible components of
$\xy$, of degrees $d_1, \dots, d_r$. Then $h$ is defined as the
composition
\[
\pi _1 (\ppe - \xy) \to H_1(\ppe - \xy, \zet)=\frac{\zet
^r}{(d_1,\dots ,d_r)} \stackrel{\alpha}{\to} \zet /d\zet ,
\]
where $\alpha $ is defined by $\alpha [(a_1,\dots ,a_r)]=\ov{\sum a_i}$.
The composition of $h$ with the characters $\rho _s:\zet /d\zet \to \cpx
^{\ast}$ defined by $\rho_s (1)=e^{2 \pi i s /d}$ (for $1\le s \le d-1$)
provide one dimensional flat bundles ${\bf V}_s$ over $\pe - \xy$ with
monodromy representation $\rho _s \circ h$.

Let $j:\ppe - \xy \hookrightarrow \ppe$ denote the inclusion map. It is
not difficult to see that the direct image sheaf $\ls _s$ ($=R^0\ls _s$)
coincides with the extension by zero $j_{!}{\bf V}_s$, $s=0,\dots
d-1$. We define the ``equivariant defect'' by:
\[
\beta_s=
\left\{
\begin{array}{l}
p_n(\xy)=b_{n-1}(\ppe - \xy) \ \ \mbox{if}\ s=0;  \\ \ \\
b_{n+1}(\ppe, \ls _s)  \ \ \ \ \ \ \ \ \ \mbox{if}\ s=1,\dots,d-1.
\end{array} \right.
\]
Here, $b_{n+1}(\ppe, \ls _s)$ is the dimension of the sheaf
cohomology
$H^{n+1}(\ppe, \ls _s)$. This vector space is the $e^{2 \pi i s/d}$-
eigenspace of $H^{n+1}(X'_0)$, where $X'_0$ is the $d$-th cyclic
covering of $\pe$ branched along $\xy$ and the action is induced by the
natural Galois action (cf. (2.12), see also \S 2, VIII).

Set $\mbox{Sing}(\xy) = \{p_1,...,p_k
\}$, and let $F_i, \mu _i, T_i$ be respectively the
local Milnor fiber, the Milnor number and the local algebraic monodromy
$H^{n}(F_i,\cpx)\to H^{n}(F_i,\cpx)$ of the isolated hypersurface
singularity $(\xy, p_i)$. We will call an invariant {\it local} if
it can be expressed in terms of the local operators $\{T_i\}_{i=1}^{k}$
and the numbers $n$ and $d$. We define the following local numerical
invariants:
\begin{eqnarray*}
\chi _0 &=& - \Sigma _{i=1}^{k}\mu_{i} + \frac{(-1)^n+(d-1)^{n+1}}{d}+
(-1)^{n+1}    \\
\chi _s &=& \chi _0 + (-1)^{n} \ \ \mbox{for} \ \ s=1,\dots , d-1.
\end{eqnarray*}

Now we are ready to formulate our main result. If $T$ is an operator,
let $T_{\alpha}$ denote its restriction to its generalized
$\alpha$-eigenspace and let $\# _l T_{\alpha}$ be the number of Jordan
blocks of $T_{\alpha}$ of size $l$. Set $\# T_{\alpha}=\Sigma_{l\ge 1}\#
_l T_{\alpha}$. With this notations one has:
\vsss

\noindent {\bf Main Theorem:}

\vspace{2mm}

\noindent {\it
I. If $\alpha=e^{2 \pi i s /d}$, $s=0,\dots,d-1$, then:
\begin{enumerate}
\item $\# _1 (\mi)_{\alpha}=\chi _s+ 2\beta_s - \Sigma_{i=1}^{k}\#
(T_i)_{\alpha}$.
\item $\# _2 (\mi)_{\alpha}= -\beta_s + \Sigma_{i=1}^{k}\# _1
(T_i)_{\alpha}$.
\item $\# _{l+1} (\mi)_{\alpha}=  \Sigma_{i=1}^{k}\# _l
(T_i)_{\alpha}$\ \ for $l\ge 2$.
\end{enumerate}}

\noindent II. {\it If $\alpha^d \neq 1$, then
$(\mi)_{\alpha}=\oplus_{i=1}^{k}\alpha^d \cdot (T_i)_{\alpha^{1-d}}$,
i.e. $\# _l (\mi)_{\alpha }= \Sigma_{i=1}^{k}\#_l(T_i)_{\alpha^{1-d}}$
for all $l\ge 1.$  }
\vsss

\noindent {\bf Corollary 1:}{\it

\begin{enumerate}
\item $\# _l (\mi)_{\alpha}=0$ for $l\ge n+2$.
\item $\# _{n+1} (\mi)_1=0$ if $\alpha^d\not=1$ or $\alpha=1$.
\end{enumerate}  }
\vsss

\noindent {\bf Corollary 2:\ }\cite[(3.3)]{GN}{\it
The characteristic polynomial of $\mi$ is given by the following local
formula:
\[
det(\lambda\cdot Id - \mi)=(\lambda - 1)^{(-1)^{n+1}} \cdot (\lambda ^d
-1)^{\frac{(d-1)^{n+1}+(-1)^n}{d}} \cdot \prod ^k_{i=1}
\frac{\mbox{det}(\lambda^{d-1}\cdot Id - T_i)}{(\lambda ^d
-1)^{\mu_{i}}}.
\]}

\noindent Another byproduct of the main theorem is the following:
\vsss

\noindent {\bf Corollary 3:\ }{\it
If $\alpha=e^{2 \pi i s /d}$, $s=0,\dots,d-1$, then:
\[
\frac{\sum_{i=1}^{k}\#(T_i)_{\alpha}-\chi _s}{2} \le \beta _s \le
\sum _{i=1}^{k} \# _1 (T_i)_{\alpha}.
\]}

If $\sum _{i=1}^{k}\#(T_i)_{\alpha}\le \chi _s$, then the lower bound
given by Corollary 3 is useless, but in some cases it gives even the
right value of $\beta_0$. For example, if $n=2$ and
$l_1,\dots,l_d\in \cpx[X,Y,Z]$
are linear forms defining an arrangement of lines in ${\bf P}^2$ such
that no two of them meet at a point, then $\beta_0=d-1$, and this is
exactly the bound given by corollary 3 above.
 From the discussion in \cite[p.161]{Di.book} follows that if $(\xy ,
p_i)$ are all non-degenerate singularities (i.e., $\Sigma
_{i=1}^k\#(T_i)_1=0$), then the defect $\beta_0=0$. Also, in
Zariski's book \cite{Zar} we can find similar criteria for $n=2$.
Notice that our bound gives a sharper criterion: If $\sum
_{i=1}^k\# _1(T_i)_1=0$ then $\beta_0=0$ (cf. also Remark 2.30 below).

For another corollary of the main theorem, see (2.16).
\vsss

\noindent {\bf \S 2. Proof of the main theorem}
\vsss

\noindent {\bf I. The main construction and two exact sequences}
\vsss

Let $f:\ea \to \cpx$ be a polynomial map which satisfies the condition
$(\ast)$. By \cite[(2.6)]{GN}, we can assume that $f$ is of the form
$f_d+x_{n+1}^{d-1}$, where $f_d$ is homogeneous of degree $d$ (and no
singularity of $\xy$ is on the hyperplane
$x_{n+1}=0$). Set (cf. \cite[\S 5]{GN}):
\[
{\cal X}=\{([x],t)\in \pe \times D :
t\Bigl(f_d(x_1,\dots,x_{n+1})+x_0x_{n+1}^{d-1}\Bigr)=x_0^d\},
\]
where $D$ denotes a disk of sufficiently small radius in the complex
plane with center at the origin. Then the map $\pi:{\cal X}\to D$ given
by $\pi([x],t)=t$ induces a locally trivial ${\cal
C}^{\infty}$-fibration over $D-\{ 0\}$ with projective fibers, these are
exactly the projective closures of the fibers of $f$. Moreover, if we
denote by $T$ the algebraic monodromy $H^{n}(\ov{f^{-1}(t_0)} ) \to
H^{n}(\ov{f^{-1}(t_0)} )$ of the projective closure $\ov{f^{-1}(t_0)}$
associated with the path $s\mapsto t_0 e^{2 \pi i s}$ ($s\in[0,1]$, $\left|
t_0 \right|$ sufficiently large), then the monodromy of $\pi$ over $\partial
D$ (with its natural orientation) is exactly $T^{-1}$. Theorem (4.6) in
\cite{GN} says basically that the knowledge of $T$ is equivalent to that
of $\mi$.
\vsss

\noindent (2.1) {\bf Theorem} \cite[(4.6)]{GN}{\it
\begin{enumerate}
\item For any $\alpha\neq 1$, $(\mi)_{\alpha }=T_{\alpha }$.
\item For $\alpha=1 $ one has
\begin{enumerate}
\item $\# _1(\mi)_1 = b_n(\xy )+p_{n-1}(\xy ) - \# T _1$.
\item $\# _2 (\mi )_1 =\# _1 T_1 - b_{n}(\xy)$.
\item $\# _{l+1}(\mi )_1 =\#_l T_1 \ \ \mbox{for} \ l\ge 2$.
\end{enumerate}
\end{enumerate}}

The big disadvantage of the map $\pi $ is that its central fiber $\pi
^{-1}(0)$ is non-reduced. For this reason we consider the following
construction:
Let $D'$ be again a disc of small radius and consider $\delta :D'\to D$
given by $\delta(t)=t^d$. Then the normalization ${\cal X}'$ of ${\cal X
}\times _{\delta}D'$ can be identified with
\[
{\cal X}'=\{([x],t)\in \pe \times D' \mid f_d(x_1,\dots ,
x_{n+1})+tx_0x_{n+1}^{d-1}=x_0^{d}\}
\]
Now $\pi ':{\cal X}'\to D'$ ($\pi([x],t)=t$) induces a locally trivial
${\cal C}^{\infty}$-fibration over $D'-\{ 0\}$ with algebraic monodromy
$T^{-d}:H^{n}((\pi ')^{-1}(t_0))\to H^{n}((\pi ')^{-1}(t_0))$.

Notice that now both ${\cal X}'$ and the central fiber $X'_0=\pi
^{-1}(0)$ have only isolated singularities: $\mbox{Sing}({\cal
X}')=\mbox{Sing}(\xy)\times \{ 0\}$. In fact, the central fiber is the
$d$-fold cyclic covering of $\ppe$ branched along $\xy$, in particular
if we set $\mbox{Sing}(X'_0)=\{ p'_1,\dots ,p'_k\}$, then the
isolated singularities $(X'_0, p'_i)$ are the
$d$-th suspensions of the singularities $(\xy , p_i)_{i=1}^k$ and the
map $\pi '$ provides their smoothings. Let $F'_i$ (respectively,
$T'_i$) be the Milnor
fiber of $(X'_0, p'_i)$ (respectively, the monodromy $H^n(F'_i)\to
H^n(F'_i)$
corresponding to the smoothing given by $\pi '$), $1\le i\le k$. Then
the exact sequence of vanishing cycles is:
\[
(2.2) \hspace{1cm}0 \to H^n(X'_0)\to H^n(X'_t) \to \oplus_{i=1}^k
H^n(F'_i) \to P^{n+1}(X'_0) \to 0
\]
where $X'_t=(\pi ')^{-1}(t)$ (for some fixed $t\neq 0$) and
$P^{n+1}(X'_0)$ is the primitive cohomology $\mbox{Ker}[H^{n+1}(X'_0)
\to H^{n+1}(X'_t)]$  (for details see \cite[(10), (5.3)]{GN}).

Our second exact sequence is given by the generalized invariant cycle
theorem proved in the Appendix of \cite{GN}:
\[
(2.3) \hspace{1cm}0\to H^n(X'_0) \to \mbox{Ker}((T^{-d})_1 - Id) \to
\oplus_{i=1}^k H^{n+1}_{\{ p'_i \}} ({\cal X}') \to 0.  \label{2.3}
\]
Both sequences are exact sequences of mixed Hodge structures and there
is a natural monodromy action on them, which at the level of
$H^{n}(X'_t)$ is $T^{-d}$. The main point of the paper is the
construction of an action on these exact sequences which at the level of
$H^n(X'_t)$ is $T^{-1}$. More precisely: we would like to understand the
monodromy of $\pi $, but this map has a non-reduced central fiber, which
makes the study difficult. Then we go to the normalization of the
$d$-fold covering $\pi $, which is $\pi '$, and we lift the monodromy of
$\pi $ to the level of $\pi '$.

First, notice that $\pi ':{\cal X}'\to D'$ has a natural Galois action
of the cyclic group $\zet /d\zet$ over $\pi:{\cal X} \to D$, that is, we
have a commutative diagram:
\vs

\begin{picture}(130,110)(-100,30)
\put(0,70){$D'$}
\put(0,130){${\cal X}'$}
\put(60,40){$D$}
\put(120,70){$D'$}
\put(120,130){${\cal X}'$}
\put(60,100){${\cal X}$}
\put(16,135){\vector(1,0){100}}
\put(16,75){\vector(1,0){100}}
\put(5,123){\vector(0,-1){40}}
\put(65,93){\vector(0,-1){40}}
\put(125,123){\vector(0,-1){40}}
\put(15,126){\vector(2,-1){37}}
\put(115,66){\vector(-2,-1){37}}
\put(115,126){\vector(-2,-1){37}}
\put(15,66){\vector(2,-1){37}}
\put(60,140){$G$}
\put(10,100){$\pi '$}
\put(130,100){$\pi '$}
\put(70,80){$\pi$}
\put(-90,90){$(2.4)$}
\end{picture}

\noindent where if we set $\xi = e^{2 \pi i /d}$, the horizontal map
$D'\to D'$ is given by $t\mapsto t\xi^{-1}$ and $G$ is given by
$G([x_0:\dots:x_{n+1}],t)=([\xi x_0:\dots : x_{n+1}], t \xi ^{-1})$.

Now we lift the geometric monodromy of $\pi $ (over $D-\{ 0\}$) to the level
of $\pi ':{\cal X}'\to D'$. Fix a point $t_0\in D'-\{ 0\}$, consider
the circle $S^1_{t_0}=\{z\in D': \left| z \right| =t_0 \}$,
and take ${\cal E}'=
(\pi ')^{-1}(S^1_{t_0})$. The fibration ${\cal E}'\to S^1_{t_0}$ is
still denoted by $\pi '$, its monodromy transformation is $T^{-d}$. Take
a local trivialization over the positive arc $[t_0,t_0\xi ]$
i.e., a diffeomorphism $h$ such that
\vs

\begin{picture}(185,90)(-40,10)
\put(25,90){$[0,\frac{1}{d}] \times X'_{t_0}$}
\put(125,100){$h$}
\put(125,85){$\simeq$}
\put(10,50){$(s,x)\mapsto t_0 e^{2 \pi i s}$}
\put(165,50){$\pi'$}
\put(120,20){$S^1_{t_0}$}
\put(185,90){$(\pi')^{-1}(\mbox{arc}[t_0,t_0\xi])$}
\put(85,95){\vector(1,0){90}}
\put(67,80){\vector(1,-1){45}}
\put(185,80){\vector(-1,-1){45}}
\put(-50,50){(2.5)}
\end{picture}

Then the geometric monodromy of $\pi$ can be identified at the level of
$\pi'$ with the composition
\[
(2.6) \hspace{3cm}X'_{t_0}\stackrel{h(\frac{1}{d},\cdot
)}{\to}X'_{t_0\xi}\stackrel{G}{\to}X'_{t_0}.\hspace{4cm}
\]
This lifting construction can be extended over $D'$ as follows: Since
$X'_0=(\pi ')^{-1}(0)$ has only isolated singularities, it is possible
to construct a flow
\vs

\begin{picture}(165,80)(-80,0)
\put(0,0){$[0,1]\times D'$}
\put(0,70){$[0,1]\times {\cal X}'$}
\put(140,0){$D'$}
\put(140,70){${\cal X}'$}
\put(55,5){$\vector(1,0){80}$}
\put(55,75){$\vector(1,0){80}$}
\put(30,60){$\vector(0,-1){45}$}
\put(145,60){$\vector(0,-1){45}$}
\put(35,35){$id \times \pi '$}
\put(150,35){$\pi '$}
\put(90,10){$\varphi$}
\put(90,80){$\phi$}
\end{picture}
\vs

\noindent such that the above diagram is commutative,
$\varphi(s,t)=te^{2 \pi i
s}$, and $\phi(s,x)=x$ for any $x\in X'_0$
(see, for example, \cite{Clemens}). Now consider the composition
$G\circ \phi(\frac{1}{d},\cdot )$ over $D'$
\vs
\vs

\begin{picture}(200,80)(-70,0)
\put(0,80){${\cal X}'$}
\put(100,80){${\cal X}'$}
\put(200,80){${\cal X}'$}
\put(20,85){\vector(1,0){75}}
\put(120,85){\vector(1,0){75}}
\put(45,95){$\phi(\frac{1}{d}, \cdot)$}
\put(150,90){$G$}
\put(32,45){$\pi '$}
\put(165,45){$\pi'$}
\put(100,10){$D'$}
\put(15,75){\vector(3,-2){80}}
\put(195,75){\vector(-3,-2){80}}
\put(-70,45){(2.7)}
\end{picture}

This will be called the ``lifted geometric monodromy". In the next
subsections we will determine the isomorphisms induced by it on the
vector spaces which appear in the exact sequences (2.2) and (2.3).
Obviously, on $H^n(X'_t)$ the induced ``lifted geometric monodromy" is
exactly $T^{-1}$.

The action on the spaces $H^q(X'_0)$ can be determined as follows.
Since $\phi(s,x)=x$ for any $x\in X'_0$, the isomorphism
$\phi(\frac{1}{d},\cdot )$ restricted to $X'_0$ is the identity.
Therefore, the action on $H^q(X'_0)$ is induced by the Galois action
$G:X'_0 \to X'_0$, $G([x_0:\dots:x_{n+1}])=[\xi x_0, \dots , x_{n+1}]$.
This action will be denoted by $G^q$.
\vsss

\noindent {\bf II. The action on $\oplus _{i=1}^k H^n(F'_{i})$.}
\vsss

If $\varphi:H\to H$ is a linear map, we will denote by $c_l(\varphi
):H^{\oplus l}\to H^{\oplus l}$ the linear map defined by $c_l(\varphi
)(x_1,\dots,x_l)=(\varphi(x_l),x_1,\dots,x_{l-1})$. Then we have:
\vsss

\noindent (2.8) {\bf Theorem:}{\it
\begin{enumerate}
\item If $S(F_i)$ denotes the suspension of $F_i$ then we have a
homotopy equivalence $F'_i \sim \bigvee_{d-1} S(F_i)$, therefore an
isomorphism $H^n(F'_i)\simeq H^{n-1}(F_i)^{\oplus(d-1)}$
\item Under the isomorphism above, the ``lifted monodromy action" on
$H^{n}(F'_i)$ is $c_{d-1}(T_i)$.
\end{enumerate}}

\noindent {\it Proof:} Part a) was proved already in \cite[(5.3)]{GN},
but it follows also from
our discussion here. For b) we use a similar construction as in [loc.
cit.]. The map-germ $({\cal X}', p'_i) \stackrel{\pi '}{\to}
(D',0)$ can be identified with
\[
{\cal Y}_i := \{ (y_0,y,t) : g_i(y)+ ty_0 =y_0^d \}
\stackrel{\pi'_i}{\to} (D',0)
\]
where $y_0$ and $y=(y_1,\dots,y_n)$ are local affine coordinates
($y_i=\frac{x_i}{x_{n+1}}$, $i=0,\dots,n$), $g_i$ is the local equation
of $(\xy, p_i)\subset(\ppe, p_i)$ and $\pi '_i$ is given by $\pi
'_i(y_0,y,t)=t$. The Galois action on ${\cal Y}_i$ is
$G(y_0,y,t)=(y_0\xi ,y, t\xi^{-1})$. As in the global situation (2.5),
consider the (local) locally trivial fibrations induced by $\pi '$,
${\cal E}'_i=(\pi ')^{-1}(S^1_{t_0})\cap {\cal Y}_i \stackrel{\pi
'_i}{\to} S^1_{t_0}$, with fiber $F'_i$. Consider the local
trivializations over the positive arc $[t_0,t_0\xi]$
\vs

\begin{picture}(185,100)(-40,0)
\put(25,90){$[0,\frac{1}{d}] \times F'_i$}
\put(125,100){$h_i$}
\put(125,85){$\simeq$}
\put(10,50){$(s,x)\mapsto t_0 e^{2 \pi i s}$}
\put(165,50){$\pi'_i$}
\put(120,20){$S^1_{t_0}$}
\put(185,90){$(\pi'_i)^{-1}(\mbox{arc}[t_0,t_0\xi])$}
\put(85,95){\vector(1,0){90}}
\put(67,80){\vector(1,-1){45}}
\put(185,80){\vector(-1,-1){45}}
\end{picture}

\noindent Then the ``lifted geometric action" on $F'_i:= (\pi
'_i)^{-1}(t_0)$ is the composition
\[
(2.9) \hspace{2cm}(\pi '_i)^{-1}(t_0) \stackrel{h_i(\frac{1}{d},
\cdot)}{\to}
(\pi '_i)^{-1}(t_0\xi) \stackrel{G}{\to} (\pi '_i)^{-1}(t_0) \hspace{2cm}
\]
which will be denoted $h'_i$. We will prove that this geometric action
induces $c_{d-1}(T_i)$ at the cohomology level.\\
{\em Remark:}\ It is not difficult to see that $(h_i')^d$ is the monodromy of
$\pi_i'$. In \cite{GN} this is identified with $c_{d-1}(T_i^d)\approx
[c_{d-1}(T_i)]^d$.

As in \cite[(5.3)]{GN}, consider the isolated complete intersection
singularity given by ${\cal Y}'_i\stackrel{\varphi}{\to}D'\times D'$,
$\varphi(y_0,y,t)=(t,y_0)$. The discriminant of $\varphi $ is $\Delta
=\{ty_0=y_0^d \}$  and $\varphi $ is a locally trivial fibration over
$D'\times D' - \Delta$ with fiber $F_i$.

For $t_0e^{2 \pi i \beta}\in S^1_{t_0}$, the intersection points of the
line $\{t=t_0e^{2 \pi i \beta}\}$ with the discriminant $\{ty_0=y_0^d\}$
of $\varphi$ are
\begin{eqnarray*}
q_0(\beta)=(t_0e^{2 \pi i \beta},0) \ \mbox{ and }\ q_j(\beta ) & = &
(t_0e^{2 \pi i \beta}, e^{2\pi i \frac{j+\beta}{d-1}}\cdot
\sqrt[d-1]{t_0})
\end{eqnarray*}
for $ j=1,\dots , d-1$.
We will use the following notations:
\begin{eqnarray*}
I_j(\beta) & = & \mbox{segment }[q_0(\beta), q_j(\beta)] \ \ (\mbox{in}\
\{t_0e^{2 \pi i \beta} \}\times D'),  \\
I(\beta) &= & \bigcup_{j=1}^{d-1}I_j(\beta ),  \\
          B& =& \bigcup _{\beta \in[0,1]}I(\beta )\subset D'\times D',\\
r_j(\beta ) & = & \mbox{middlepoint of } I_j(\beta )=
            (t_0e^{2 \pi i \beta }, {1\over 2} \sqrt[d-1]{t_0}\,
            e^{2 \pi i \frac{j+\beta }{d-1}}).
\end{eqnarray*}
It is obvious that, for all $1\le j \le d-1$, $\varphi ^{-1}(q_j(\beta
))$ is contractible
and $\varphi ^{-1}(r_j(\beta ))$ is exactly $F_i$.
Therefore $\varphi^{-1}(I_j(\beta))$ can be identified with the
suspension $S(F_i)$ and $\varphi^{-1}(I(\beta)) \sim
\bigvee_{d-1}S(F_i)$.

The inclusion $B\subset S^1_{t_0} \times D'$ admits a strong deformation
retract which can be lifted. Consider the torus ${\bf T}=S^1_{t_0}\times
S^1 _{\frac{1}{2} \sqrt[d-1]{t_0}}$ which contains the points
$r_j(\beta )$. By the identification of $F'_i=(\pi
'_i)^{-1}(t_0)=\varphi^{-1}(\{ t_0 \}\times D')$ with $\varphi
^{-1}(I(0))\sim \bigvee _{d-1}S(F_i)$, the homology of $F'_i$ is
generated by a wedge of suspensions of cycles which lie above the points
$r_j(0)$, $1\le j \le d-1$. When we move $t$ on the positive arc
$[t_0,t_0\xi]$, then these points move on the path $[0,\frac{1}{d}]\to
{\bf T}$ given by $\beta \mapsto r_j (\beta )$. We denote these paths
by $\gamma_j$, with endpoints $r_j(0)$ and $r_j(\frac{1}{d})$, i.e.
\[
\gamma_j (s) = ( t_0 e^{2 \pi i s}, {1\over 2} \sqrt[d-1]{t_0}\ e^{2
\pi i \frac{j+s}{d-1}}) \ , \ s\in[0,\frac{1}{d}].
\]
The local trivialization over $\cup_j\gamma_j$ corresponds to
$h_i(\frac{1}{d}, \cdot )$ in (2.9) (we will explain this identification
more precisely later). Next, we identify the Galois action with
some local trivialization over some paths.

\noindent Consider the paths $\tau_j : [0, \frac{1}{d}] \to {\bf T}$ defined by
\[
\tau_j (s)=(t_0e^{2\pi i (\frac{1}{d}-s)}, {1\over 2} \sqrt[d-1]{t_0}\
e^{2 \pi i (\frac{j-1}{d-1}+\frac{1}{d(d-1)}+s)}),
\]
which connects $r_{j-1}(\frac{1}{d})$ and $r_{j}(0)$.
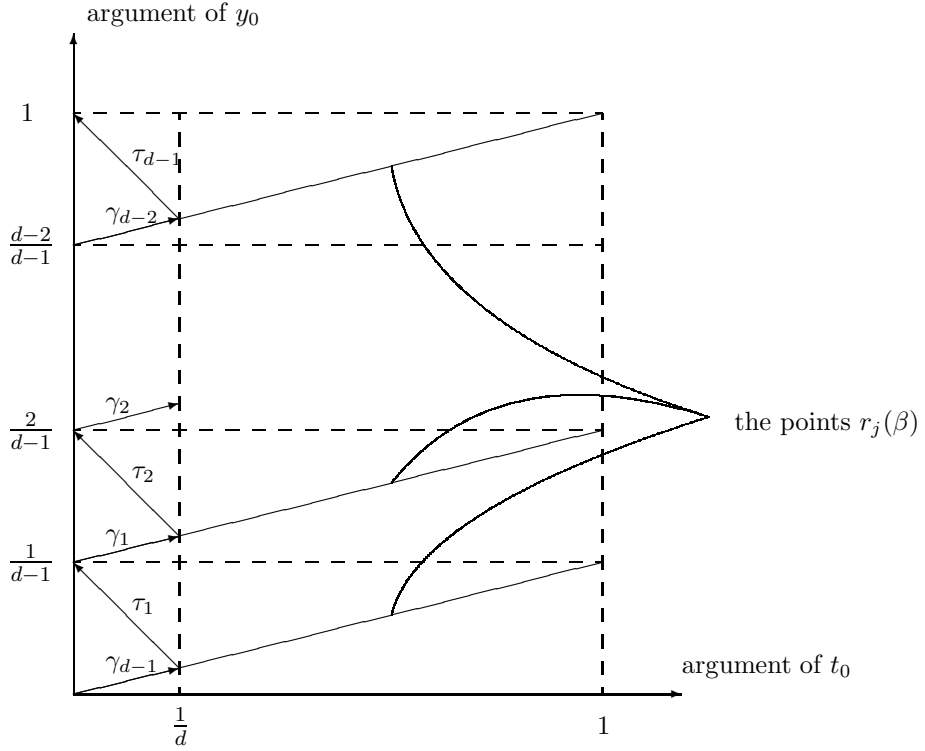
\vs

\begin{picture}(270,270)(-40,0)
\put(20,20){\vector(0,1){250}}
\put(20,20){\vector(1,0){230}}

\put(20,20){\vector(4,1){40}}
\put(60,30){\vector(-1,1){40}}
\put(20,70){\vector(4,1){40}}
\put(60,80){\vector(-1,1){40}}
\put(20,120){\vector(4,1){40}}
\put(20,190){\vector(4,1){40}}
\put(60,200){\vector(-1,1){40}}
\put(20,20){\line(4,1){200}}
\put(20,70){\line(4,1){200}}
\put(20,190){\line(4,1){200}}

\put(20,20){\dashbox{5}(200,50)}
\put(20,20){\dashbox{5}(200,100)}
\put(20,20){\dashbox{5}(200,170)}
\put(20,20){\dashbox{5}(200,220)}
\put(20,20){\dashbox{5}(40,220)}

\put(57,5){$\frac{1}{d}$}
\put(218,5){\footnotesize{$1$}}
\put(-5,67){$\frac{1}{d-1}$}
\put(-5,117){$\frac{2}{d-1}$}
\put(-5,187){$\frac{d-2}{d-1}$}
\put(0,237){\footnotesize{$1$}}

\put(250,28){\footnotesize{argument of $t_0$}}
\put(25,275){\footnotesize{argument of $y_0$}}
\put(270,120){\footnotesize{the points $r_j(\beta )$}}

\put(32,30){\footnotesize{$\gamma _{d-1}$}}
\put(32,78){\footnotesize{$\gamma _{1}$}}
\put(32,128){\footnotesize{$\gamma _{2}$}}
\put(32,200){\footnotesize{$\gamma _{d-2}$}}

\put(42,52){\footnotesize{$\tau _{1}$}}
\put(42,102){\footnotesize{$\tau _{2}$}}
\put(42,222){\footnotesize{$\tau _{d-1}$}}

 \bezier{300}(260,125)(150,90)(140,50)
 \bezier{300}(260,125)(180,150)(140,100)
 \bezier{300}(260,125)(150,160)(140,220)
\end{picture}

\noindent Notice that the Galois action $(y_0,y,t) \mapsto (y_0
\xi, y, t \xi^{-1})$ induces
\[
G:\varphi^{-1}(r_{j-1}(\frac{1}{d})) \stackrel{\sim}{\to}
\varphi^{-1}(r_{j-1}(0)).
\]
Consider the isomorphism (up to isotopy) given by the local
trivialization of $\varphi $ above the oriented path $\tau _{j}$:
\[
Tr_j:
\varphi^{-1}(r_{j-1}(\frac{1}{d})) \stackrel{\sim}{\to}
\varphi^{-1}(r_{j-1}(0)).
\]

\noindent {\it Fact:} The composition
\[
(2.10) \hspace{1cm}\varphi^{-1}(r_{j-1}(0)) \stackrel{G^{-1}}{\to}
\varphi^{-1}(r_{j-1}(\frac{1}{d})) \stackrel{Tr_j}{\to}
\varphi^{-1}(r_{j-1}(0))
\]
is isotopic to the identity.

\noindent {\it Proof of the fact:} Consider the map $\delta : D'\times
D'\to D'$
given by $\delta(t,y_0)=y_0^d-ty_0$ (one has $\delta ^{-1}(0)=$the
discriminant of $\varphi$). First notice that
$\delta(r_{j-1}(\frac{1}{d}))=\delta(r_j(0))$. Therefore, the
composition (2.10) can be identified with
\[
g_i^{-1}(\delta(r_j(0)))\stackrel{Id}{\to}
g_i^{-1}(\delta(r_{j-1}(\frac{1}{d})))\stackrel{\ti{Tr_j}}{\to}
g_i^{-1}(\delta(r_j(0)))
\]
where the first map is the identity $y\mapsto y$ (the second component
of $G$) and $\ti{Tr_j}$ is the trivialization of $g_i$ above the loop
$s\mapsto \delta(\tau_k(s))$, ($s\in[0,\frac{1}{d}]$). Now, it is easy
to verify that $\delta(\tau_k(s))$ can be written in the form $B \cdot
e^{2 \pi i (a+ds)} + A$, where $\left| A \right| > \left| B \right|$.
Therefore the loop
$s\mapsto \delta(\tau_k(s))$ is isotopic to zero in $D'- \{0\}$. So,
$\ti{Tr}_j$ (and hence $Tr_j \circ G^{-1}$ too) is isotopic to the
identity. $\Box$
\vs

The above fact shows that the Galois action $G$ can be replaced by the
local trivialization above the paths $\{\tau_j\} _j$.

Since $h(\frac{1}{d}, \cdot )$ in (2.9) corresponds to the local
trivialization above $\{ \gamma _j\}_j$ and the Galois action to the
local trivialization above $\{\tau_j\} _j$, then the composed map in
(2.9) corresponds to the trivialization above $\{\tau_{j+1}\circ \gamma
_j\}_j$. Now we identify the fibers of $\varphi $ above the points
\[
r_{d-1}(0),r_{d-1}(\frac{1}{d}),r_1(0),r_1(\frac{1}{d}), \dots ,
r_{d-1}(0),r_{d-2}(\frac{1}{d})
\]
via the paths:
\[
(2.11) \hspace{2cm}\gamma _{d-1}, \tau _1 , \gamma _1 , \tau _2 , \gamma
_2 , \tau _3, \dots , \gamma _{d-2} . \hspace{2cm}
\]
The fiber $F'_i$ is
\[
S\varphi^{-1}(r_{d-1}(0)) \vee S\varphi^{-1}(r_{1}(0)) \vee \dots \vee
S\varphi^{-1}(r_{d-2}(0))
\]
and the ``lifted monodromy action" is induced by
\[
S(\tau_1\circ\gamma_{d-1})_{\ast}\vee
S(\tau_2\circ\gamma_1)_{\ast}\vee
\dots \vee S(\tau_{d-1}\circ\gamma_{d-2})_{\ast}
\]
But this  (because of the identification of fibers via the
paths in (2.11)) is exactly the isomorphism $c_{d-1}(Q)$, where $Q$ is the
monodromy of $\varphi$ above the loop
\[
l=\gamma _{d-1}\circ\tau_1 \circ \gamma_{1}\circ\tau_2 \circ \dots \circ
\gamma_{d-2}\circ \tau_{d-1}.
\]
The loop $l$ in the complement of the discriminant of $\varphi $ is
homotopic to $s\mapsto (t_0,e^{2 \pi i s}\frac{\sqrt[d-1]{t_0}}{2})$,
$s\in[0,1]$. The linking number of this (second) loop with $\{y_0=0\}$
is one, and with $\{y_0^{d-1}=t_0\}$ is zero. Therefore $Q=T_i$, in
particular the map $h'_{i}$ induced by $G\circ h_i(\frac{1}{d},
\cdot)$ is $c_{d-1}(T_i)$.  $\Box$
\vs

\noindent {\bf III. The action on $H^{n+1}_{\{p'_i\}}({\cal X}')$}.
\vs

In this subsection we prove that the ``lifted action" on
$H^{n+1}_{\{p'_i\}}({\cal X}')$ is trivial. Let ${\cal K}'_i$ be the
link of $({\cal X}',p'_i)=\{g_i(y)+ty_0-y_0^{d}=0\}$ (we use the same
notations as in II). The map
\[
\pi':\{ g_i(y)+ ty_0 - y_0^{d}=0 \} \to D', \ \ \
(y_0,y,t)\mapsto t
\]
gives an open book decomposition of ${\cal K}'_i$. Let
$K_i=\{ t=0 \} \subseteq {\cal K}'_i$ be the link of $t$, then
$arg=arg(t): {\cal K}'_i - K_i \to S^1$ is a ${\cal C}^{\infty}$-locally
trivial fibration. Consider the flow $\phi :[0,1]\times {\cal K}'_i \to
{\cal K}'_i$ such that
\begin{enumerate}
\item If $x\in\{ t=0 \}=K_i$, then $\phi(s,x)=x$ for any $s$.
\item If $x\not\in K_i$, then $arg(\phi(s,x))=e^{2 \pi i s} arg(s)$.
\end{enumerate}
The wanted geometric action is the composed map
\[
{\cal K}'_i \stackrel{\phi(\frac{1}{d}, \cdot )}{\to} {\cal K}'_i
\stackrel{G}{\to} {\cal K}'_i
\]
Now $\phi (\frac{1}{d}, \cdot )$ is isotopic to the identity via the
flow
$\phi(s,\cdot)$, $s\in[0,\frac{1}{d}]$. The Galois action is isotopic to
the identity as well. To see this consider the isotopy:
\begin{eqnarray*}
(s,(y_0,y,t))&\mapsto &(y_0(s),y(s),t(s))= \\
& = & (y_0e^{2 \pi i s /d},y,(t-y_0^{d-1})e^{-2 \pi i s
/d}+y_0^{d-1}e^{2 \pi i s (d-1)/d}).
\end{eqnarray*}
If $s=0$, then $(y_0(0),y(0),t(0))=(y_0,y,t)$, if $s=1$ then
$(y_0(0),y(0),t(0))=(y_0\xi,y,t\xi^{-1})=G(y_0,y,t)$.
\vs

\noindent {\bf IV. The exact sequences revisited.}
\vs

We summarize the results of the subsections I-III: One has the following
two exact sequences, with the ``lifted monodromy action":
\vs

\begin{picture}(370,65)(0,0)
\put(0,60){$0$}
\put(40,60){$H^n(X'_0)$}
\put(115,60){$H^n(X'_t)$}
\put(190,60){$\oplus_{i=1}^kH^n(F'_i)$}
\put(285,60){$P^{n+1}(X'_0)$}
\put(370,60){$0$}
\put(10,65){\vector(1,0){25}}
\put(85,65){\vector(1,0){25}}
\put(160,65){\vector(1,0){25}}
\put(255,65){\vector(1,0){25}}
\put(340,65){\vector(1,0){25}}

\put(0,10){$0$}
\put(40,10){$H^n(X'_0)$}
\put(115,10){$H^n(X'_t)$}
\put(190,10){$\oplus_{i=1}^kH^n(F'_i)$}
\put(285,10){$P^{n+1}(X'_0)$}
\put(370,10){$0$}
\put(10,15){\vector(1,0){25}}
\put(85,15){\vector(1,0){25}}
\put(160,15){\vector(1,0){25}}
\put(255,15){\vector(1,0){25}}
\put(340,15){\vector(1,0){25}}

\put(60,50){\vector(0,-1){25}}
\put(135,50){\vector(0,-1){25}}
\put(220,50){\vector(0,-1){25}}
\put(310,50){\vector(0,-1){25}}

\put(225,35){\small{${\oplus \ c_{d-1}(T_i)}$}}
\put(65,35){$G^n$}
\put(140,35){$T^{-1}$}
\put(315,35){$G^{n+1}$}
\put(-20,35){(E.1)}
\end{picture}

and (E.2):
\vs

\begin{picture}(370,65)(0,0)

\put(0,60){$0$}
\put(40,60){$H^n(X'_0)$}
\put(115,60){$\mbox{Ker}\ (T^{-d}-Id)$}
\put(230,60){$\oplus_{i=1}^k H^{n+1}_{\{p'_i\}}({\cal X}')$}
\put(340,60){$0$}

\put(10,65){\vector(1,0){25}}
\put(85,65){\vector(1,0){25}}
\put(200,65){\vector(1,0){25}}
\put(310,65){\vector(1,0){25}}

\put(0,10){$0$}
\put(40,10){$H^n(X'_0)$}
\put(115,10){$\mbox{Ker}\ (T^{-d}-Id)$}
\put(230,10){$\oplus_{i=1}^k H^{n+1}_{\{p'_i\}}({\cal X}')$}
\put(340,10){$0$}

\put(10,15){\vector(1,0){25}}
\put(85,15){\vector(1,0){25}}
\put(200,15){\vector(1,0){25}}
\put(310,15){\vector(1,0){25}}

\put(60,50){\vector(0,-1){25}}
\put(155,50){\vector(0,-1){25}}
\put(260,50){\vector(0,-1){25}}

\put(65,35){$G^n$}
\put(160,35){$T^{-1}$}
\put(265,35){\small{Identity}}
\end{picture}

The main theorem will follow from these exact sequences and from
some mixed Hodge-theoretical arguments.

We end this subsection with some facts about the Galois action
$G^{\ast}:H^{\ast}(X'_0) \to H^{\ast}(X'_0)$, where
$p:X'_0 \to \ppe$ is the $d$-th cyclic covering branched along $\xy$.
One has that $H^q(X'_0,\cpx)={\bf H}^q(\ppe ,{\bf R}p_{\ast}\cpx
_{X'_0})=
H^q(\ppe , p_{\ast}\cpx _{X'_0})$ and the restriction of $p_{\ast}\cpx
_{X'_0}$ to the complement of $\xy $ is a flat bundle. Its corresponding
monodromy representation is given by the composed map $\pi _1 (\ppe -
\xy)\stackrel{h}{\to} \zet /d\zet \stackrel{r}{\to} \mbox{Aut}(\zet ^d)$
(see \S 1 for the definition of $h$), where $r(1):= \sigma
:\zet/d\zet\to\zet/d\zet$ is the permutation $\sigma (x_1,\dots
,x_d)=(x_d, x_1, \dots ,x_{d-1})$.

But also the Galois action is induced by $\sigma $. So, we have a direct sum
decomposition $p_{\ast }\cpx_{X'_0} = \cpx _{\ppe}\oplus
\oplus_{s=1}^{d-1}\ls _{s}$ such that $G\bigm|_{\cpx_{\ppe}}$ is the
identity and $G\bigm|_{\ls _s}$ is the multiplication by $\xi^s=e^{2
\pi i s /d}$. Therefore:
\[
(2.12) \hspace{1cm}(H^q(X'_0); G^q)=(H^q(\ppe)\oplus
(\oplus_{s=1}^{d-1}H^q(\ppe, \ls _s)); \oplus_{s=0}^{d-1}\xi^s).
\]

\vspace{3mm}

\noindent {\bf V. The proof of the main Theorem, case $\alpha=1$.}
\vs

Consider the exact sequence (E.2) with its actions. Since the
(generalized) $1$-eigenspace of $T^{-1}$ on $\mbox{Ker}(T^{-d}-Id)$ is
exactly $\mbox{Ker}(T^{-1}-Id)$, the decomposition (2.12) provides the
following exact sequence:
\[
(2.13) \hspace{1cm}0 \to H^n(\ppe) \to \mbox{Ker}(T^{-1}-Id) \to
\oplus_{i=1}^{k}H^{n+1}_{\{p'_i\}}({\cal X}') \to 0
\]
This is again an exact sequence of mixed Hodge structures. Now, it is
on the one hand clear that the weight of $H^n(\ppe)$ is $n$ and, on the
other hand,
\[
\mbox{dim}\ Gr^W_{n-l+1}H^{n+1}_{\{p'_i\}}({\cal X}')=\#_l(T_i)_1  \ \
(\mbox{for } l\in\zet)
\]
(see \cite[(5.5)]{GN}). Since the weight filtration on $H^n(X'_t)$ is
the monodromy weight filtration of $T^{-1}$ centered at $n$, one has
$\mbox{dim}\ Gr^W_{n-l+1}\mbox{Ker}(T^{-1}-Id)=\#_l (T^{-1})$. This
shows that
\[
(2.14) \hspace{1cm}(T^{-1})_1 = \oplus_{i=1}^k (T_i)_1 \oplus
\left\{
\begin{array}{l}
0 \ \ \ \ \ \ \mbox{ if } \  $n$ \ \mbox{is odd}  \\
id_{\cpx} \ \ \  \mbox{ if } \ $n$ \ \mbox{is even}
\end{array} \right.
\]
Now notice that
\[
(2.15) \hspace{2cm}p_{n-1}(\xy) - p_{n}(\xy) = \chi_0 \hspace{2cm}
\]
(see, for example (2.29), or \cite{Di.book}).
Hence the result follows from (2.1), (2.14) and (2.15).
$\Box$

\noindent (2.16) {\it Remark:}
Using the exact sequence of vanishing cycles of $X^{\infty}\subset {\bf P}^n$,
by standard mixed Hodge theoretical arguments, one can prove that the
dimensions $\dim Gr^W_{n-l}P^{n-1}(X^{\infty})$ \ $(l\geq 1)$  are equal to the
numbers  on the right hand side of the equalities in the Main Theorem,
case $\alpha=1$.
Hence , the main theorem gives:
$$\#_l(T^{\infty}_f)_1=\dim Gr^W_{n-l}P^{n-1}(X^{\infty})\ \ \mbox{for} \ \
l\in {\bf Z},$$
where $W$ denotes the weight filtration.
\vs

\noindent {\bf VI. The proof of the main Theorem, case $\alpha^d=1,
\alpha\neq 1$.}
\vs

First notice that the Galois action $G:X'_0\to X'_0$ is an algebraic
map, therefore $G^q:H^q(X'_0)\to H^q(X'_0)$ preserves the weight
filtration.  Since $\cls{n}$ is the $\xi^{s}$-eigenspace of $G^q$ (cf.
2.12), it has a natural induced weight filtration (actually, it has a
natural mixed Hodge structure). Now let $\alpha=e^{2 \pi i s /d}=\xi^s$
for $1\le s\le d-1$. Then by (E.2) one has:
\[
(2.17) \hspace{1cm}\cls{n}=\mbox{Ker}(T^{-1}-\alpha).
\]
In particular, for $l\ge 1$ one has
\[
(2.18)
\hspace{1cm}\#_l(T^{-1})_{\alpha}=\mbox{dim}\ Gr^W_{n-l+1}\cls{n}.
\]
Actually, (2.18) together with (2.1) already determine
$(\mi)_{\alpha}$, but this is not exactly the assertion of the main
theorem, we want some more information about the right hand side of
(2.18).

Consider the exact sequence (E.1). Using (2.17) one has:
\vs

\begin{picture}(370,65)(15,0)
\put(0,60){$0$}
\put(35,60){Ker$(T^{-1}-\alpha)$}
\put(135,60){$H^n(X'_t)_{\alpha}$}
\put(205,60){$[\oplus H^n(F'_i)]_{\alpha}$}
\put(290,60){$H^{n+1}(\ppe,\ls _s)$}
\put(400,60){$0$}
\put(10,65){\vector(1,0){20}}
\put(110,65){\vector(1,0){20}}
\put(180,65){\vector(1,0){20}}
\put(265,65){\vector(1,0){20}}
\put(375,65){\vector(1,0){20}}

\put(0,10){$0$}
\put(35,10){Ker$(T^{-1}-\alpha)$}
\put(135,10){$H^n(X'_t)_{\alpha}$}
\put(205,10){$[\oplus H^n(F'_i)]_{\alpha}$}
\put(290,10){$H^{n+1}(\ppe,\ls _s)$}
\put(400,10){$0$}

\put(10,15){\vector(1,0){20}}
\put(110,15){\vector(1,0){20}}
\put(180,15){\vector(1,0){20}}
\put(265,15){\vector(1,0){20}}
\put(375,15){\vector(1,0){20}}

\put(70,50){\vector(0,-1){25}}
\put(155,50){\vector(0,-1){25}}
\put(230,50){\vector(0,-1){25}}
\put(330,50){\vector(0,-1){25}}

\put(75,35){$\cdot \alpha$}
\put(160,35){$(T^{-1})_{\alpha}$}
\put(235,35){$[\oplus c_{d-1}(T_i)]_{\alpha}$}
\put(335,35){$\cdot \alpha$}
\end{picture}

By \cite{Sch}, the weight filtration of $H^n(X'_t)_{\alpha}$ is the
monodromy weight
filtration of $(T^{-1})_{\alpha}$ centered at $n$, thus the quotient
$H^n(X'_t)/\mbox{Ker}\ (T^{-1}-\alpha )$ has a (polarized) mixed Hodge
structure with weight filtration equal to the monodromy weight
filtration of the class $[T^{-1}_{\alpha}]$ of $(T^{-1})_{\alpha}$
centered at $n+1$. On the other hand, $H^{n+1}(X'_0)$ is pure of weight
$n+1$ (\cite{Stee.Oslo}) and $G^{n+1}$ preserves the weight filtration,
hence $\cls{n+1}$ is pure of weight $n+1$ (cf. (2.12)). These two facts
show that the weight filtration of $[\oplus H^n(F'_i)]_{\alpha}$ is the
monodromy weight filtration of $(\oplus c_{d-1}(T_i))_{\alpha }$. Now
comparing the dimensions of the primitive cohomologies of
$[T^{-1}_{\alpha}]$ and $(c_{d-1}(T_i))_{\alpha}$ one has:
\[
(2.19) \hspace{1cm}\left\{
\begin{array}{l}
\#_2 (T^{-1})_{\alpha} = - \mbox{dim}\ H^{n+1}(\ppe, \ls _s) +
\Sigma_{i=1}^k \#_1(c_{d-1}(T_i))_{\alpha} \\ \ \\
\#_{l+1}(T^{-1})_{\alpha}=\Sigma_{i=1}^k \#_l(c_{d-1}(T_i))_{\alpha} \
\mbox{ for }\ l\ge 2.
\end{array} \right.
\]
Since $\alpha^d=1$ and $\alpha\neq 1$,
$(c_{d-1}(T_i))_{\alpha}=(T_i)_{\alpha^{-1}}$. Therefore:
\[
(2.20) \hspace{1cm}\left\{
\begin{array}{l}
\#_2(T^{-1})_{\alpha}=-\beta_s + \Sigma _{i=1}^k \#_1
(T_i)_{\alpha^{-1}} \\ \ \\
\#_{l+1}(T^{-1})_{\alpha}=\Sigma_{i=1}^k \#_{l}(T_i)_{\alpha^{-1}} \ \mbox{ for
}\ l\ge 2.
\end{array} \right.
\]
Since $\beta _{d-s}=\beta_s$ \ ($s=1,\dots ,d-1$), parts (Ib,Ic) follow
from (2.20) and (2.1.a). In order to prove (Ia) (i.e., to compute
$\#_1(\mi)_{\alpha }=\#_1 T_{\alpha}$) notice that from (2.18) one has:
\[
(2.21) \hspace{1cm}\#(T^{-1})_{\alpha}=\mbox{dim}\ H^n(\ppe,\ls _s),
\]
and from (2.20):
\[
(2.22) \hspace{1cm}\Sigma_{l\ge 2}\#_l(T^{-1})_{\alpha }=-\beta _s +
\Sigma _{i=1}^k \# (T_i)_{\alpha ^{-1}}.
\]
Therefore
\[
\#_1(T^{-1})_{\alpha}= \mbox{dim}\ H^n(\ppe, \ls _s) + \mbox{dim}\
H^{n+1}(\ppe, \ls _s) - \Sigma_{i=1}^k \# (T_i)_{\alpha ^{-1}}.
\]
Now (Ia) follows from this identity and the following result:
\vs

\noindent (2.23) {\bf Proposition:} {\it For $s=1,\dots ,d-1$ one has:
\begin{enumerate}
\item dim $H^q(\ppe, \ls _s)=0$ \ \ for $q\neq n,n+1$.
\item $\mbox{dim}\ H^n(\ppe, \ls _s)-\mbox{dim}\ H^{n+1}(\ppe, \ls
_s)=\chi _s$.
\end{enumerate}
(For the definition of $\chi _s$, see the introduction). In
particular, the Euler characteristic of $(\ppe, \ls _s)$ is local.
}
\vs

\noindent {\it Proof:} The first part follows from (2.12), because
$P^q(X'_0)=0$
if $q\neq n,n+1$. The second part will be proved (together with some
other relations) in subsection VIII, (2.29).
\vs

%
\noindent {\bf VII. The proof of the main theorem, case $\alpha^d\neq
1$.} \vs

Consider the exact sequence (E.1). Since $(G^{\ast})^d=Id$, the
generalized $\alpha $-eigenspaces are:
\[
(2.25) \hspace{1cm}(H^n(X'_t)_{\alpha}, (T^{-1})_{\alpha}) \simeq
((\oplus_{i=1}^k H^n(F'_i))_{\alpha }, (\oplus c_{d-1}(T_i))_{\alpha}).
\]
By (2.1), $T_{\alpha}=(\mi)_{\alpha}$. Also, if $\Pi_j(\lambda - \xi_j)$
is the characteristic polynomial of $T_i$, then $\Pi_j(\lambda ^{d-1} -
\xi_j)$ is the characteristic polynomial of $c_{d-1}(T_i)$. Therefore,
if $\alpha $ is an eigenvalue of $c_{d-1}(T_i)$, then $\alpha^{d-1}=\xi
_j$ for some $\xi_j $ and the unipotent (or nilpotent) part of
$c_{d-1}(T_i)_{\alpha}$ and $(T_i)_{\xi_j}$ can be identified. This ends
the proof of the main theorem.
\vs

\noindent {\bf VIII. The relation with the Milnor fiber of $f_d:\ea \to
\cpx$.} \vs

Consider the homogeneous singularity $f_d:\ea \to \cpx$ with
one-dimensional singular locus, let $F$ be its Milnor fiber. It is
well-known that its (reduced) homology is concentrated in $H_n(F)$ and
$H_{n-1}(F)$. Let $h_q:H_q(F)\to H_q(F)$ be the algebraic monodromy of
$f_d$, where $q=n-1,n$. In this subsection we identify our local and global
invariants with numerical invariants given by the transformations
$h_{n-1}, h_n$.

Recall that we denoted $p:X'_0\to\ppe$ the $d$-th cyclic covering
branched along $\xy$.
Then $p^{-1}(\xy)$ can be identified with $\xy$ and $X'_0 -p^{-1}(\xy)$
with $F$. This gives a cyclic (unramified) covering \ $F\to\ppe -\xy$
with
fiber $\zet/d\zet$. By duality, $H_q(F)=H^{2n-q}(X'_0,p^{-1}(\xy))$,
therefore the exact sequence of the pair $(X'_0,p^{-1}(\xy))$ reads as
follows:
\vs

\begin{picture}(370,65)(20,0)
\put(0,60){$0$}
\put(30,60){$P^{n-1}(\xy)$}
\put(110,60){$H_n(F)$}
\put(170,60){$P^n(X'_0)$}
\put(235,60){$P^n(\xy)$}
\put(305,60){$H_{n-1}(F)$}
\put(370,60){$P^{n+1}(X'_0)\to 0$}

\put(10,65){\vector(1,0){15}}
\put(90,65){\vector(1,0){15}}
\put(150,65){\vector(1,0){15}}
\put(215,65){\vector(1,0){15}}
\put(285,65){\vector(1,0){15}}
\put(355,65){\vector(1,0){15}}

\put(0,10){$0$}
\put(30,10){$P^{n-1}(\xy)$}
\put(110,10){$H_n(F)$}
\put(170,10){$P^n(X'_0)$}
\put(235,10){$P^n(\xy)$}
\put(305,10){$H_{n-1}(F)$}
\put(370,10){$P^{n+1}(X'_0)\to 0$}

\put(10,15){\vector(1,0){15}}
\put(90,15){\vector(1,0){15}}
\put(150,15){\vector(1,0){15}}
\put(215,15){\vector(1,0){15}}
\put(285,15){\vector(1,0){15}}
\put(355,15){\vector(1,0){15}}

\put(60,50){\vector(0,-1){25}}
\put(130,50){\vector(0,-1){25}}
\put(190,50){\vector(0,-1){25}}
\put(260,50){\vector(0,-1){25}}
\put(330,50){\vector(0,-1){25}}
\put(400,50){\vector(0,-1){25}}

\put(65,35){$Id$}
\put(135,35){$h_n^{-1}$}
\put(195,35){$G^n$}
\put(265,35){$Id$}
\put(335,35){$h^{-1}_{n-1}$}
\put(405,35){$G^{n+1}$}
\put(-20,35){(2.26)}
\end{picture}

In the above diagram, we have also inserted the corresponding Galois
actions. The Galois action on $X'_0=\{f_d(x_1, \dots ,x_{n+1})=x_0^d \}$
is $[x_0:\dots :x_{n+1}]\mapsto [\xi x_0:\dots :x_{n+1}]$. If on $X'_0 -
p^{-1}(\xy)$ we take affine coordinates $y_i=x_i/x_0\ , \ 1\le i \le
n+1$, then the induced action is $(y_1,\dots ,y_{n+1})\mapsto
\xi^{-1}(y_1,\dots ,y_{n+1})$. This is the inverse of the geometric
monodromy of the Milnor fiber $F$.

Now, if we consider the generalized eigenspaces of the Galois action in
(2.26), one has the following identifications:
\[
(2.27) \hspace{1cm}\left\{
\begin{array}{l}
(H_n(F), h_n^{-1})_{\neq 1} = (P^n(X'_0), G^n)_{\neq 1}  \\ \ \\
(H_{n-1}(F), h_{n-1}^{-1})_{\neq 1} = (P^{n+1}(X'_0), G^{n+1})_{\neq 1}
\\ \ \\
\mbox{dim}\ H_{n}(F)_1 = p_{n-1}(\xy) \ \mbox{and \
dim}\ H_{n-1}(F)_1=p_n(\xy).
\end{array} \right.
\]
We recall (cf. (2.12)) that $(P^q(X'_0),G^q)_{\neq
1}=(\oplus_{s=1}^{d-1}H^q(\ppe,\ls _s), \oplus_{s=1}^{d-1}\xi^s)$. In
particular, all our global invariants are equivalent to the
characteristic polynomial of $h_{n-1}$, i.e.
 $\beta_s=\mbox{rank} H_{n-1}(F)_{\alpha}$ ($0\leq s\leq d-1$).
Now let us consider
the zeta function of $f_d:\ea \to \cpx$. This is basically given in
\cite{Siersma2}:
\begin{eqnarray*}
(2.28) \hspace{.4cm}\frac{\mbox{det}(\lambda \cdot Id -
h_n)}{\mbox{det}(\lambda \cdot Id
- h_{n-1})}& = &(\lambda -1)^{(-1)^{n+1}} \cdot (\lambda ^d -
1)^{\frac{(d-1)^{n+1}+(-1)^{n}}{d}}\cdot \prod_{i=1}^k(\lambda^d -
1)^{-\mu_{i}} \\
& = & \prod_{s=0}^{d-1}(\lambda - e^{2 \pi i s/d})^{\chi_s}.
\end{eqnarray*}
Now, (2.28) and (2.27) give:
\[
(2.29) \hspace{.3cm}\left\{
\begin{array}{l}
p_{n-1}(\xy) - p_n(\xy)=\chi_0 \\ \ \\
\mbox{dim}\ H^n(\ppe, \ls _s) - \mbox{dim}\ H^{n+1}(\ppe, \ls _s)=\chi_s
\ \ (s=1,\dots ,d-1).
\end{array} \right.
\]
This proves (2.15) and (2.23.b).
\vs

\noindent (2.30){\it Remark:} By (2.27) and corollary 3 one has (for
$s=0,\dots ,d-1$) \[
\mbox{dim}\ H_{n-1}(F)_{e^{2\pi i s/d}}=\beta _s\le
\Sigma_{i=1}^k\#_1(T_i)_{e^{2 \pi i s/d}}
\]
Similar restrictions can be found in \cite[(5.4) and \S 9]{Siersma}.
\vs

\noindent (2.31)\ {\em Remark:}\
{\em  (The relation with $\pi_{n-1}({\bf P}^n-X^{\infty})$.)}\\
Here  we present  the connection between the present paper and
\cite{L}, more precisely, between  the defects $\beta_s\ (0\leq s\leq d-1)$
and $\pi_{n-1}({\bf P}^n-X^{\infty})$.

First, assume that $n=2$. Denote $G=\pi_1({\bf P}^n-X^{\infty}),\
G'=[G,G]$, and $G''=[G',G']$.
Then $0\to \pi_1(F)\to G\to {\bf Z}/d{\bf Z}\to 0$
is an exact sequence, actually $\pi_1(F)=G'$.
Therefore, $H_1(F)=G'/G''$, and it
has a natural action of ${\bf Z}/d{\bf Z}$. By (2.27) one has:
$$\beta_s=\mbox{rank}((G'/G'')\otimes {\bf Q})_{\alpha};\ \ \alpha=e^{2\pi
is/d};
\ s=0,\ldots, d-1.$$
Now, assume that $n>2$. From the covering $F\to {\bf P}^n-X^{\infty}$ and
Hurewicz theorem one has: $\pi_1({\bf P}^n-X^{\infty})={\bf Z}/d{\bf Z}$,
$\pi_q({\bf P}^n-X^{\infty})=0$ if $1<q<n-1$, and $\pi_{n-1}({\bf
P}^n-X^{\infty})=\pi_{n-1}(F)=H_{n-1}(F)$. Hence by (2.27):
$$\beta_s=\mbox{rank} (\pi_{n-1}({\bf P}^n-X^{\infty})
\otimes {\bf Q})_{\alpha};\ \ \alpha=e^{2\pi is/d};
\ s=0,\ldots, d-1,$$
where the action of $\pi_1({\bf P}^n-X^{\infty})$ on
$\pi_{n-1}({\bf P}^n-X^{\infty})$ is the natural one.

We thank Professor A. Libgober bringing in our attention  the invariant
$\pi_{n-1}$, and his helpful comments about it.

\vs

\noindent {\bf \S 3. Examples}
\vs

\noindent {\it I. Zariski's plane sextics:}
\vs

Set $d=6$ and let
$f_6\in\cpx[X;Y;Z]$ be a form defining a plane
sextic in ${\bf P}^2$ with six cusps and no other singularities. Then
$\chi_0=8$ and $\chi_s=9$ for $s=1,\dots ,5$. Moreover (since the
characteristic polynomial of the local monodromy of a cusp singularity
is $t^2-t+1$), $\beta_s=0$ if $s=0,2,3,4$. Our main theorem gives:
\begin{enumerate}
\item If $\alpha^d\neq 1$, then $(\mi)_{\alpha}$ has only
one-dimensional Jordan blocks, and $(\mi)_{\alpha}=Id_{\cpx^6}$ if
$\alpha=e^{\pi i \varphi}$, $\varphi\in\{\frac{1}{15}, \frac{7}{15},
\frac{11}{15}, \frac{13}{15}, \frac{17}{15}, \frac{19}{15},
\frac{23}{15}, \frac{29}{15}\}$ (i.e., if $\alpha^5=e^{2 \pi i s /6}$
for $s=1$ or $s=5$ and $\alpha^6\neq 1$). Otherwise $(\mi)_{\alpha}=0$.
\item $(\mi)_{e^{2\pi i s /6}}$ has only one-dimensional blocks if
$s=0,2,3,4$. The number of them is $8$ if $s=0$ and $9$ if $s=2,3,4$.
\item $(\mi)_{e^{2 \pi i s/6}}$ has only one and two dimensional blocks
if $s=1$ or $s=5$. The number of one-dimensional blocks is $3+2\beta_s$,
and the number of two dimensional blocks is $6-\beta_s$.
\end{enumerate}

\noindent Now, by the identification (2.27) and \cite[Theorem
2.9]{Di.book} one
has $\beta_s=1, \ (s=1,5)$ if the cusps are on a conic and $\beta_s=0$
otherwise.

\vs

\noindent {\it II. Nodal hypersurfaces:}
\vs

Assume that $f_d$ defines a hypersurface in $\ppe $ with only
nodal (i.e. $A_1$) singularities, let $k$ denote the number of nodes.
It follows from the main theorem that the maximal size of a Jordan block
of $\mi$ is two. The numbers $\beta_s$ can be computed using
\cite[VI, Theorem (4.5)]{Di.book} and (2.27).

If $dn$ is even,
set $S=\pr$, \ $q=\frac{dn}{2}-n-1$, and let $S_q$ denote the
homogeneous
component of degree $q$ of $S$. If $\Sigma \subset \xy$ denotes the set
of nodes of $\{f_d=0\}$,
let $S_q(\Sigma)=\{h\in S_q \mid h_{\mid \Sigma}=0 \}$, and
$\mbox{defect} (S_q(\Sigma)) := k - \mbox{codim}_{S_q}(S_q(\Sigma)).$

{}From the main
theorem we get the following possibilities for $\mi$: \begin{enumerate}
\item $n$ is odd, $d$ is odd.

Here $\beta_s=0$ for all
$s$. Thus $\mi$ has no Jordan blocks of size two, i.e. it is of finite
order.
\item $n$ is odd, $d$ is even.

In this case $\beta_s=0$ for $s\neq \frac{d}{2}$ and
$\mi$ can have Jordan blocks of size two only for eigenvalue
$-1$, the number of them is $\#_2(\mi)_{-1}=k-\beta _{d/2}
=k-\mbox{defect}(S_q(\Sigma))=
\mbox{codim}_{S_q}(S_q(\Sigma))$.
\item $n$ is even.

In this case, $\beta_s=0$ for $s\neq
0$ and $\beta_0=p_n(\xy)$. It follows that $\mi$
can have Jordan blocks of size two only for eigenvalue $1$, and
$\#_2 (\mi)_1=k-p_n(\xy)
=k-\mbox{defect}(S_q(\Sigma))=
\mbox{codim}_{S_q}(S_q(\Sigma))$.
This number, in general, is not zero. For example,  the defect$(S_5(\Sigma))
=\beta_0$
of the recently constructed  quintic hypersurface in ${\bf P}^4$ with $k=130$
nodes \cite{Duco} is 29 [loc.cit., p. 864].
The quintic constructed by Hirzebruch \cite{H} has 126 nodes and defect
$\beta_0=25$. Actually, there are quintics in ${\bf P}^4$ with
118 nodes and defect $18\leq \beta_0\leq 19$ \cite{W}.

\end{enumerate}

\begin{footnotesize}

\end{footnotesize}
\end{document}